\documentclass{amsart}
\usepackage{sty}

\usepackage{tikz}
\usetikzlibrary{shapes}
\usetikzlibrary{decorations.markings}
\tikzset{every picture/.style={>=latex}}
\tikzstyle{vert}=[scale=.8,
      style={circle,fill=black,minimum size=5pt,inner sep=0pt}]
\tikzstyle{close}=[inner sep=0,minimum size=0]
\tikzstyle{label}=[scale=.9]
\tikzstyle{thinedge}=[->,color=green]
\tikzstyle{thinedgerev}=[<-,color=green]
\tikzstyle{thickedge}=[->,very thick]
\tikzstyle{thickedgerev}=[<-,very thick]
\tikzset{variable/.default=} 
\tikzset{label/.default=} 

\usepackage{graphicx}
\usepackage{float}
\usepackage{hyperref}\hypersetup{colorlinks,citecolor=blue,linktocpage}

\newcommand{\simp}{\operatorname{Simp}}

\usepackage{comment}
\excludecomment{longversion}
\includecomment{showfigures}

\author[Valery Alexeev,
Ryan Livingston, 
Joseph Tenini,
et al.]%
{Valery Alexeev,
Ryan Livingston,
Joseph Tenini,
Maxim Arap,
Xiaoyan Hu,
Lauren Huckaba,
Patrick Mcfaddin,
Stacy Musgrave,
Jaeho Shin,
and Catherine Ulrich}

\address{Department of Mathematics, University of Georgia}

\title[Extended Torelli map]%
{Extended Torelli map to the\\Igusa blowup in genus 6, 7, and 8}
\date{May 14, 2011}

\begin{document}
\maketitle

\begin{abstract}
  It was conjectured in \cite{Namikawa_ExtendedTorelli} that the
  Torelli map $\M_g\to\A_g$ associating to a curve its jacobian
  extends to a regular map from the Deligne-Mumford moduli space of
  stable curves $\oM_g$ to the (normalization of the) Igusa blowup
  $\oA_g\cent$. A counterexample in genus $g=9$ was found in
  \cite{AlexeevBrunyate}. Here, we prove that the extended map is
  regular for all $g\le8$, thus completely solving the problem in
  every genus.
\end{abstract}


\section{Introduction}

The Torelli map $\M_g \to \A_g$ associates to a smooth curve $C$ its
jacobian $JC$, a principally polarized abelian variety.  Does it
extend to a regular map $\oM_g\to \oA_g$, where $\oM_g$ is
Deligne-Mumford's moduli space of stable curves, and $\oA_g$ is a
toroidal compactification of $\A_g$?

This question was first asked in a pioneering paper of Namikawa
\cite{Namikawa_ExtendedTorelli} in the case when $\oA_g=\oA_g^{\rm cent}$ is
the normalization of the Igusa blowup $\Bl_{\partial A_g^*}A_g^*$ of the
Satake compactification along the boundary. This compactification was
introduced by Igusa in \cite{Igusa_BlowUp}, and is possibly the first
toroidal compactification ever constructed. It corresponds to the
central cone decomposition.

Namikawa conjectured that the extended map is always regular. He was
able to prove it for the stable curves with a planar dual graph, and
for curves of low genus; the case $g\le6$ was stated without
proof. (Note: the graphs in this paper may have multiple edges and loops).

The question was recently revisited in \cite{AlexeevBrunyate}, who
showed the following:

\medskip\noindent (1) Let $C$ be a stable curve of genus $g$, and let
$\Gamma$ be its dual graph. Then the rational map $\oM_g\to\oA_g^{\rm
  cent}$ is regular in a neighborhood of the point $[C]\in\oM_g$
$\iff$ there exists a positive definite integral-valued quadratic form
$q$ on the first cohomology $H^1(\Gamma, \bZ)$ such that $q(e_i^*)=1$
for every non-bridge edge $e_i$ of~$\Gamma$.  Such quadratic forms $q$
are called \emph{integral edge-minimizing metrics} or \emph{$\bZ$-emms} for
short.

Recall that for a graph, $H^1(\Gamma) = C^1(\Gamma) / dC^0(\Gamma)$,
where $C^1(\Gamma,\bZ)=\oplus_{{\rm edges}\ e_i}\bZ e_i^*$,
$C^0(\Gamma,\bZ)=\oplus_{{\rm vertices}\ v_j}\bZ v_j^*$, and $dv_j^* =
\sum_{e_i {\rm\ begins\ with\ } v_j } e_i^* - 
\sum_{e_i {\rm\ ends\ with\ } v_j } e_i^*$. 
We denote the image of $e_i^*$ in $H^1(\Gamma,\bZ)$ by the same letter
$e_i^*$ and call it a \emph{coedge}.

\medskip\noindent
(2) Call a graph \emph{cohomology-irreducible} if there does not exist a
partition of its edges into two groups $I_1\sqcup I_2$ such that
\begin{math}
  H^1(\Gamma,\bZ) = 
  \langle e^*_i,\ i\in I_1 \rangle
  \oplus
  \langle e^*_i,\ i\in I_2 \rangle.
\end{math}
Then $\Gamma$ is either a simple loop (a graph with one vertex and one
edge), or $\Gamma$ is loopless and 2-connected.

For every graph $\Gamma$, one has $H^1(\Gamma,\bZ) = \oplus
H^1(\Gamma_k, \bZ)$ for some cohomology-irreducible graphs $\Gamma_k$
and all coedges $e_i^*$ lie in the direct summands.
We call $G_k$ \emph{cohomology-irreducible components of}
$\Gamma$.  Then there exists a $\bZ$-emm for $\Gamma$ $\iff$ there
exist $\bZ$-emms for all $\Gamma_k$.

\medskip\noindent 
(3) If a graph $\Gamma$ is cohomology-irreducible and $q$ is a
$\bZ$-emm for $\Gamma$ then the lattice $(H^1(\Gamma,\bZ),2q)$ is a
root lattice of type $A_g$, $D_g$ ($g\ge4$), or $E_g$ ($g=6,7,8$).
Further, there exists a $\bZ$-emm of type $A_g$ $\iff$ $\Gamma$ is
planar, and for $g\ge4$ there exists a $\bZ$-emm of type $D_g$ $\iff$
$\Gamma$ is projective planar, i.e. can be embedded into the
projective plane $P=\mathbb R\mathbb P^2$.

\medskip The famous theorem of Kuratowski says that a graph is
non-planar iff it contains a subgraph homeomorphic either to $K_5$ or
to $K_{3,3}$. A Kuratowski-type theorem for the projective plane $P$
was proved by Archdeacon \cite{Archdeacon_ProjPlane,Archdeacon_Thesis}
who showed that the list of 103 minimal non-projective planar graphs
produced earlier by Glover-Huneke-Wang
\cite{GloverHunekeWang_103graphs} is complete; any other
non-projective planar graph contains a subgraph homeomorphic to one of
them. The smallest graph on their list has genus~6.

This implies that every graph of genus $\le 5$ has a $\bZ$-emm, and
consequently the extended Torelli map $\oM_g\to\oA_g^{\rm cent}$ is
regular for $g\le5$. On the other hand, as \cite{AlexeevBrunyate}
notes, there exist cohomology-irreducible non-projective planar graphs 
of genus 9, so the extended Torelli map is \emph{not} regular for every
$g\ge9$.

\medskip
Here are the main results of this paper:

\begin{theorem}\label{thm:main}
  Let $\Gamma$ be a cohomology-irreducible non-projectively planar graph of
  genus $g=6, 7$, or 8. Then $\Gamma$ admits a $\bZ$-emm of type
  $E_g$.
\end{theorem}

\begin{corollary}
  The extended Torelli map $\oM_g\to\oA_g^{\rm cent}$ is regular for
  $g\le 8$.
\end{corollary}

\begin{corollary}
  Let $C$ be a stable curve of genus $g$, and $\Gamma$ be its dual
  graph. Then the extended Torelli map $\oM_g\to\oA_g^{\rm cent}$ is
  regular in a neighborhood of the point $[C]\in\oM_g$ $\iff$ every
  cohomology-irreducible component $\Gamma_k$ has genus $\le8$ or is a
  projectively planar graph of genus $\ge9$.
\end{corollary}

The plan of the paper is as follows. In Section~\ref{sec:reduction},
we reduce the proof of Theorem~\ref{thm:main} to checking finitely
many graphs: one graph for $g=6$, 14 graphs for $g=7$, and 2394 graphs
for $g=8$. In Section~\ref{sec:6}, we give a finite algorithm for an
arbitrary graph, and then run it for the only graph needed in genus
6. In Section~\ref{sec:7}, we give the 14 graphs in genus 7 that have to
be checked, and explicitly list a $\bZ$-emm for each of them. In
Section~\ref{sec:8}, we state our computer-aided findings for genus 8.

\begin{acknowledgments}
  The subject of the paper was one of the topics of a VIGRE research
  group at the University of Georgia in the Fall of 2010, led by the
  first author. We would like to acknowledge NSF's VIGRE support under
  DMS-0738586. We thank Boris Alexeev for writing an alternative code
  for finding $\bZ$-emms for a graph.
\end{acknowledgments}

\section{Reduction to finitely many graphs}
\label{sec:reduction}

As noted in \cite[Sec.2]{AlexeevBrunyate}, for the proof of
Theorem~\ref{thm:main} we may reduce to graphs which are trivalent. So
let $H$ be a cohomology-irreducible non-projectively planar 
trivalent graph of genus $g=6,7$
or $8$. One says that $H$ is \emph{irreducible with respect to $P$} if $H$
does not embed into $P$, but for any edge $e$ in $H$, $H-e$ does embed
into $P$.  We now describe a process which will reduce $H$ to a
trivalent graph irreducible w.r.t. $P$. The operations (3a), (3b), (3c) are
illustrated in Figures~\ref{fig-a}, \ref{fig-b}, \ref{fig-c}.

\begin{enumerate}
\item If the graph is irreducible w.r.t. $P$, stop and call
  this graph $H'$.
\item If not, choose an edge $e$ so that $H-e$ does not embed 
into $P$ and delete $e$ from the graph.
  \item[(3a)] If $e$ was not a loop and did not have a parallel edge, 
then, denoting by $v_1$ and $v_2$ the distinct vertices to which $e$ 
is incident, contract an edge incident to $v_1$ and an edge incident 
to $v_2$. 
  \item[(3b)] If $e$ was not a loop but had a parallel edge $f$, then,
denoting by $v_1$ and $v_2$ the distinct vertices to which $e$ and $f$ 
are incident, contract the edge incident to $v_1$ and different from $f$
and the edge incident to $v_2$ and different from $f$.
  \item[(3c)] If $e$ was a loop incident to $v$, then delete the 
remaining edge $f$ incident to $v$ and, denoting by $w$ the other 
vertex to which $f$ is incident, contract one of the other two edges 
incident to $w$ and different from $f$.
\end{enumerate}

Notice that the above operations (3a)-(3b)-(3c) drop the genus
of the graph by 1 except for operation (3a) when e is a bridge.
Repeating this process we get a graph
$H'$ irreducible w.r.t. $P$ which is of the form 
$H'=\wH \cup \{u_1, \dots ,u_k\}$ where
the $u_i$ are isolated vertices and $\wH$ is a trivalent graph 
irreducible w.r.t. $P$.  By
\cite{GloverHuneke_CubicIrrGraphs,Milgram_IrrGraphs2} (see also
\cite{Archdeacon_ProjPlane,Archdeacon_Thesis}), $\wH$ is
isomorphic to one of the following:
\begin{enumerate}
\item[(i)] The connected graph $G$ of genus 6 shown in Figure~\ref{fig-g}.
\item[(ii)] The connected graphs $F_{11},F_{12},F_{13},F_{14}$ of genus 7
  shown in Figures~\ref{fig-f11}-\ref{fig-f14}.
\item[(iii)] The graph $E_{42}$ shown in Figure~\ref{fig-e42}.
\end{enumerate}

\begin{showfigures}
\begin{figure}[H]
 \begin{tikzpicture}
  [scale=.7,auto=left]
  \node[vert] (n1) at (1,4) {};
  \node[vert] (n2) at (3,6)  {};
  \node[vert] (n3) at (5,4)  {};
  \node[vert] (n4) at (3,2) {};
  \node[vert] (n5) at (3,4)  {};
  \node[vert] (n6) at (7,7)  {};
  \node[vert] (n7) at (11,4) {};
  \node[vert] (n8) at (13,6) {};
  \node[vert] (n9) at (15,4) {};
  \node[vert] (n10) at (13,2) {};
  \node[vert] (n11) at (13,4) {};
  \node[vert] (n12) at (9,7) {};

  \draw[-] (n1) to [out=90,in=180] (n2);
  \draw[-] (n2) to [out=0,in=90] (n3);
  \draw[-] (n3) to [out=-90,in=0] (n4);
  \draw[-] (n4) to [out=180,in=-90] (n1);

  \draw[-] (n7) to [out=90,in=180] (n8);
  \draw[-] (n8) to [out=0,in=90] (n9);
  \draw[-] (n9) to [out=-90,in=0] (n10);
  \draw[-] (n10) to [out=180,in=-90] (n7);

  \draw[-] (n1) -- (n5);
  \draw[-] (n5) -- (n3);
  \draw[-] (n7) -- (n11);
  \draw[-] (n11) -- (n9);

  \draw[-] (n2) to [out=45,in=180] (n6);
  \draw[-] (n5) to [out=60,in=210] (n6);
  \draw[-] (n4) to [out=-30,in=-90] (n6);
  \draw[-] (n8) to [out=165,in=0] (n12);
  \draw[-] (n11) to [out=120,in=-30] (n12);
  \draw[-] (n10) to [out=210,in=-90] (n12);

\end{tikzpicture}
\caption{The Graph $E_{42}$.}
\label{fig-e42}
\end{figure}
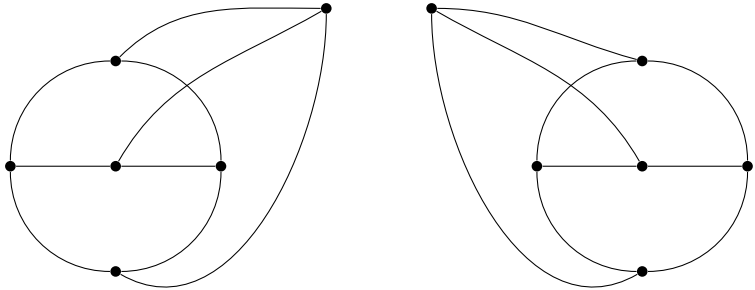
\end{showfigures}

Thus, we may construct $H$ from
$\wH$ by reversing the algorithm above. We make this explicit for
the relevant genera $6,7$ and $8$.

\subsection{$H$ has genus 6}
Since $H$ is cohomology-irreducible, it has no bridges and so operations
(3a), (3b) and (3c) would all drop the genus. Thus $H$ is already 
irredicible with respect to $P$ and so $H=\wH=G$.  Thus, to show the existence
of $\bZ$-emms for graphs of genus~$6$, it suffices to produce one for
$G$.

\subsection{$H$ has genus 7}
Either $\wH$ equals one of
$F_{11},F_{12},F_{13},F_{14}$ or $\wH=G$. In the first case we
have that $H$ is equal to one of $F_{11},F_{12},F_{13},F_{14}$ (again since $H$ 
was cohomology-irreducible, thus bridgeless).  The
second case is slightly more complicated. First notice that $H'$ has
at most one isolated vertex $v$, because in the case of applying (3c), 
the genus drops by 1. Then $H$ may be obtained from $\wH$ by doing one of the
following three operations. Notice that (a), (b) and (c) are the 
inverse operations of (3a), (3b) and (3c) (defined above) 
respectively. 

\begin{enumerate}
\item[(a)] Choose two distinct edges $e_1$ and $e_2$ and add an edge from 
the midpoint of $e_1$ to the midpoint of $e_2$.
\item[(b)] Choose an edge and add a handle to it.
\item[(c)] Choose an edge $e'$ and add an edge $f$ from the midpoint of $e'$ to
the isolated vertex $v$. Then add a loop $e$ to $v$.
\end{enumerate}

\begin{showfigures}
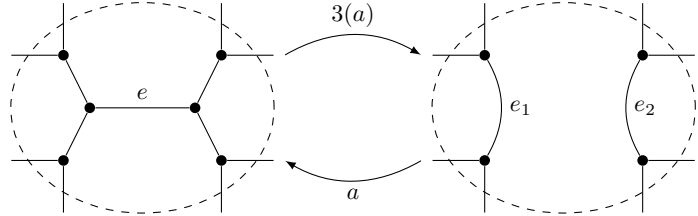
\begin{figure}[H]
 \begin{tikzpicture}
  [scale=.7,auto=left]

  \draw[dashed,thin] (3.5,3) ellipse (2.5 and 2);
  \draw[dashed,thin] (11.5,3) ellipse (2.5 and 2);

  \node[close] (n1) at (1,2) {};
  \node[close] (n2) at (2,1)  {};
  \node[vert] (n3) at (2,2)  {};
  \node[vert] (n4) at (2.5,3) {};
  \node[vert] (n5) at (2,4)  {};
  \node[close] (n6) at (1,4)  {};
  \node[close] (n7) at (2,5) {};
  \node[vert] (n8) at (4.5,3) {};
  \node[vert] (n9) at (5,2) {};
  \node[close] (n10) at (5,1) {};
  \node[close] (n11) at (6,2) {};
  \node[vert] (n12) at (5,4) {};
  \node[close] (n13) at (5,5) {};
  \node[close] (n14) at (6,4) {};
  
  \node[close] (n15) at (6.2,4) {};
  \node[close] (n16) at (8.8,4) {};

  \node[close] (n29) at (6.2,2) {};
  \node[close] (n30) at (8.8,2) {};

  \node[close] (n17) at (9,2) {};
  \node[close] (n18) at (10,1)  {};
  \node[vert] (n19) at (10,2)  {};
  \node[vert] (n20) at (10,4)  {};
  \node[close] (n21) at (9,4)  {};
  \node[close] (n22) at (10,5) {};
  \node[vert] (n23) at (13,2) {};
  \node[close] (n24) at (13,1) {};
  \node[close] (n25) at (14,2) {};
  \node[vert] (n26) at (13,4) {};
  \node[close] (n27) at (13,5) {};
  \node[close] (n28) at (14,4) {};

  \draw[-] (n1) -- (n3);
  \draw[-] (n2) -- (n3);
  \draw[-] (n3) -- (n4);
  \draw[-] (n4) -- (n5);
  \draw[-] (n5) -- (n6);
  \draw[-] (n5) -- (n7);
  \draw[-] (n4) -- node[label]{$e$} (n8);
  \draw[-] (n8) -- (n9);
  \draw[-] (n9) -- (n10);
  \draw[-] (n9) -- (n11);
  \draw[-] (n8) -- (n12);
  \draw[-] (n12) -- (n13);
  \draw[-] (n12) -- (n14);

  \draw[->] (n15) to [out=30,in=150] node[label]{$3(a)$} (n16);
  \draw[->] (n30) to [out=210,in=-30] node[label]{$a$} (n29);

  \draw[-] (n17) -- (n19);
  \draw[-] (n18) -- (n19);
  \draw[-] (n20) -- (n21);
  \draw[-] (n20) -- (n22);
  \draw[-] (n23) -- (n24);
  \draw[-] (n23) -- (n25);
  \draw[-] (n26) -- (n27);
  \draw[-] (n26) -- (n28);

  \draw[-] (n20) to [out=-60,in=60] node[label]{$e_1$} (n19);
  \draw[-] (n26) to [out=-120,in=120] node[label]{$e_2$} (n23);

\end{tikzpicture}
\caption{The procedures (3a) and (a).}
\label{fig-a}
\end{figure}
\end{showfigures}

\begin{showfigures}
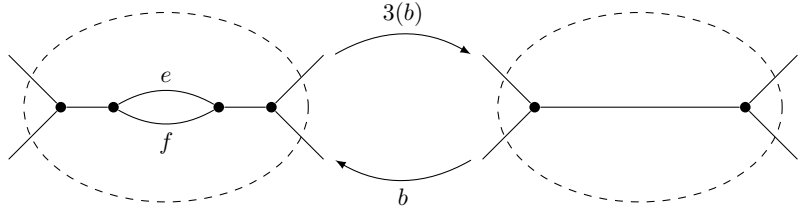
\begin{figure}[H]
 \begin{tikzpicture}
  [scale=.7,auto=left]

  \draw[dashed,thin] (3,3) ellipse (2.7 and 1.8);
  \draw[dashed,thin] (12,3) ellipse (2.7 and 1.8);

  \node[close] (n1) at (0,4) {};
  \node[close] (n2) at (0,2) {};
  \node[vert] (n3) at (1,3) {};
  \node[vert] (n4) at (2,3) {};
  \node[vert] (n5) at (4,3) {};
  \node[vert] (n6) at (5,3) {};
  \node[close] (n7) at (6,4) {};
  \node[close] (n8) at (6,2) {};

  \node[close] (n9) at (9,4) {};
  \node[close] (n10) at (9,2) {};
  \node[vert] (n11) at (10,3) {};
  \node[vert] (n12) at (14,3) {};
  \node[close] (n13) at (15,4) {};
  \node[close] (n14) at (15,2) {};

  \node[close] (n15) at (6.2,4) {};
  \node[close] (n16) at (8.8,4) {};

  \node[close] (n17) at (6.2,2) {};
  \node[close] (n18) at (8.8,2) {};

  \draw[-] (n1) -- (n3);
  \draw[-] (n2) -- (n3);
  \draw[-] (n3) -- (n4);
  \draw[-] (n4) to [out=30,in=150] node[label]{$e$} (n5);
  \draw[-] (n4) to [out=-30,in=210] node[label][swap]{$f$} (n5);
  \draw[-] (n5) -- (n6);
  \draw[-] (n6) -- (n7);
  \draw[-] (n6) -- (n8);

  \draw[->] (n15) to [out=30,in=150] node[label]{$3(b)$} (n16);
  \draw[->] (n18) to [out=210,in=-30] node[label]{$b$} (n17);

  \draw[-] (n9) -- (n11);
  \draw[-] (n10) -- (n11);
  \draw[-] (n11) -- (n12);
  \draw[-] (n12) -- (n13);
  \draw[-] (n12) -- (n14);

\end{tikzpicture}
\caption{The procedures (3b) and (b).}
\label{fig-b}
\end{figure}
\end{showfigures}

\begin{showfigures}
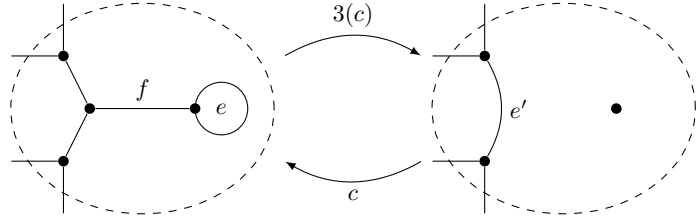
\begin{figure}[H]
 \begin{tikzpicture}
  [scale=.7,auto=left]

  \draw[dashed,thin] (3.5,3) ellipse (2.5 and 2);
  \draw[dashed,thin] (11.5,3) ellipse (2.5 and 2);
  \draw (5,3) circle (.5) node[label]{$e$};

  \node[close] (n1) at (1,2) {};
  \node[close] (n2) at (2,1)  {};
  \node[vert] (n3) at (2,2)  {};
  \node[vert] (n4) at (2.5,3) {};
  \node[vert] (n5) at (2,4)  {};
  \node[close] (n6) at (1,4)  {};
  \node[close] (n7) at (2,5) {};
  \node[vert] (n8) at (4.5,3) {};
  
  \node[close] (n15) at (6.2,4) {};
  \node[close] (n16) at (8.8,4) {};

  \node[close] (n29) at (6.2,2) {};
  \node[close] (n30) at (8.8,2) {};

  \node[close] (n17) at (9,2) {};
  \node[close] (n18) at (10,1)  {};
  \node[vert] (n19) at (10,2)  {};
  \node[vert] (n20) at (10,4)  {};
  \node[close] (n21) at (9,4)  {};
  \node[close] (n22) at (10,5) {};
  \node[vert] (n23) at (12.5,3) {};
  
  \draw[-] (n1) -- (n3);
  \draw[-] (n2) -- (n3);
  \draw[-] (n3) -- (n4);
  \draw[-] (n4) -- (n5);
  \draw[-] (n5) -- (n6);
  \draw[-] (n5) -- (n7);
  \draw[-] (n4) -- node[label]{$f$} (n8);

  \draw[->] (n15) to [out=30,in=150] node[label]{$3(c)$} (n16);
  \draw[->] (n30) to [out=210,in=-30] node[label]{$c$} (n29);

  \draw[-] (n17) -- (n19);
  \draw[-] (n18) -- (n19);
  \draw[-] (n20) -- (n21);
  \draw[-] (n20) -- (n22);

  \draw[-] (n20) to [out=-60,in=60] node[label]{$e'$}(n19);

\end{tikzpicture}
\caption{The procedures (3c) and (c).}
\label{fig-c}
\end{figure}
\end{showfigures}

In the case (c), $f$ is a bridge and so we do not need to consider 
graphs acquired from $\wH$ from operation (c).

A careful but elementary analysis shows that the cases (a) and (b), up
to symmetries, produce ten possible graphs for $H$. We denote these
graphs $G_1, \dots G_{10}$, they appear in figures
\ref{fig-g1}-\ref{fig-g10} below.  Thus, to show the existence of
$\bZ$-emms for graphs of genus $7$ it suffices to produce one for
$F_{11},F_{12},F_{13},F_{14}$ and $G_{i}$ for $i \in \{1, \dots
,10\}$.

\begin{longversion} 

In case it is needed, below is the `careful but elementary analysis' which 
produces $G_1$ through $G_{10}$.

Below we have taken the graph $G$ and created an equivalence relation 
on the edge set. Two edges $e$ and $e'$ are equivalent if there is an 
automorphism $\phi$ of $G$ such that 
$\phi(e)=e'$.

\begin{figure}[H]
\begin{tikzpicture}
 [scale=.7,auto=left]
 \node[vert] (n1) at (1,4) {};
 \node[vert] (n2) at (3,6)  {};
 \node[vert] (n3) at (5,4)  {};
 \node[vert] (n4) at (3,2) {};
 \node[vert] (n5) at (3,4)  {};
 \node[vert] (n7) at (8,4) {};
\node[vert] (n8) at (10,6) {};
 \node[vert] (n9) at (12,4) {};
 \node[vert] (n10) at (10,2) {};
 \node[vert] (n11) at (10,4) {};

 \draw[-](n1) to [out=90,in=180]  node[label]{$[2]$} (n2);
 \draw[-] (n2) to [out=0,in=90]  node[label]{$[2]$} (n3);
 \draw[-] (n3) to [out=-90,in=0]  node[label]{$[2]$} (n4);
 \draw[-] (n4) to [out=180,in=-90]  node[label]{$[2]$} (n1);

 \draw[-] (n7) to [out=90,in=180]  node[label]{$[2]$} (n8);
 \draw[-] (n8) to [out=0,in=90]  node[label]{$[2]$} (n9);
 \draw[-] (n9) to [out=-90,in=0]  node[label]{$[2]$} (n10);
 \draw[-] (n10) to [out=180,in=-90]  node[label]{$[2]$} (n7);

 \draw[-] (n1) --  node[label]{$[2]$} (n5);
 \draw[-] (n5) --  node[label][swap]{$[2]$} (n3);
 \draw[-] (n7) --  node[label][swap]{$[2]$} (n11);
 \draw[-] (n11) -- node[label]{$[2]$} (n9);

 \draw[-] (n2) to [out=30,in=150] node[label]{$[1]$} (n8);
 \draw[-] (n5) to [out=30,in=150] node[label]{$[1]$} (n11);
 \draw[-] (n4) to [out=-30,in=-150] node[label]{$[1]$} (n10);

\end{tikzpicture}
\caption{$G$ with edges labeled by equivalence class.}
\end{figure}

Thus, we can get two non-isomorphic graphs from $G$ by adding a handle: Let 
$G_{9}$ the graph obtained by adding a handle to edge $[1]$ and let 
$G_{10}$ be the graph obtained by adding a handle to edge $[2]$.

The final step will be to list the isomorphism classes of graphs obtained from $G$
by adding an edge joining the midpoints of two distinct edges. For this we proceed 
as follows: We must select an edge, create a vertex at its midpoint and relabel the 
equivalence classes of edges of this new graph (we will not label the edges incident
to the midpoint). Selecting an edge in $[1]$ we get 
$\wG_{[1]}$ and selecting an edge in $[2]$ we get $\wG_{[2]}$. Now for 
each equivalence class $[i]$ of edges in $\wG_{[j]}$ we can obtain a new graph 
by adding a vertex in the middle of an edge in class $[i]$ and adding an edge incident 
to the two midpoint vertices.

\begin{figure}[H]
\begin{tikzpicture}
 [scale=.7,auto=left]
 \node[vert] (n1) at (1,4) {};
 \node[vert] (n2) at (3,6)  {};
 \node[vert] (n3) at (5,4)  {};
 \node[vert] (n4) at (3,2) {};
 \node[vert] (n5) at (3,4)  {};
 \node[vert] (n7) at (9,4) {};
 \node[vert] (n8) at (11,6) {};
 \node[vert] (n9) at (13,4) {};
 \node[vert] (n10) at (11,2) {};
 \node[vert] (n11) at (11,4) {};
 \node[vert] (n12) at (7,1) {};

 \draw[-] (n1) to [out=90,in=180] node[label]{$[2]$} (n2);
 \draw[-] (n2) to [out=0,in=90] node[label]{$[2]$}  (n3);
 \draw[-] (n3) to [out=-90,in=0] node[label]{$[3]$} (n4);
 \draw[-] (n4) to [out=180,in=-90] node[label]{$[3]$} (n1);

 \draw[-] (n7) to [out=90,in=180] node[label]{$[2]$} (n8);
 \draw[-] (n8) to [out=0,in=90] node[label]{$[2]$} (n9);
 \draw[-] (n9) to [out=-90,in=0] node[label]{$[3]$} (n10);
 \draw[-] (n10) to [out=180,in=-90] node[label]{$[3]$} (n7);

 \draw[-] (n1) -- node[label]{$[2]$} (n5);
 \draw[-] (n5) -- node[label][swap]{$[2]$} (n3);
 \draw[-] (n7) -- node[label][swap]{$[2]$} (n11);
 \draw[-] (n11) -- node[label]{$[2]$} (n9);

 \draw[-] (n2) to [out=30,in=150] node[label]{$[1]$} (n8);
 \draw[-] (n5) to [out=30,in=150] node[label]{$[1]$} (n11);
 \draw[-] (n4) to [out=-30,in=180] node[label][swap]{} (n12);
 \draw[-] (n12) to [out=0,in=-160] node[label][swap]{} (n10);

\end{tikzpicture}
\caption{The graph $\wG_{[1]}$}
\end{figure}

Since we have 3 equivalence classes of edges in $\wG_{[1]}$, we may
obtain 3 graphs from $\wG_{[1]}$. The graph obtained
by choosing $[1]$ we call $G_1$. The graph obtained by choosing $[2]$
we call $G_8$. Finally, the graph obtained by choosing $[3]$ we call
$G_7$.

Suppose now that we choose an edge in the class $[2]$, we get the
graph $\wG_{[2]}$ below, where edges have been relabeled according to
their new equivalence classes.

\begin{figure}[H]
\begin{tikzpicture}
 [scale=.7,auto=left]
 \node[vert] (n1) at (1,4) {};
 \node[vert] (n2) at (3,6)  {};
 \node[vert] (n3) at (5,4)  {};
 \node[vert] (n4) at (3,2) {};
 \node[vert] (n5) at (3,4)  {};
 \node[vert] (n7) at (8,4) {};
 \node[vert] (n8) at (10,6) {};
 \node[vert] (n9) at (12,4) {};
 \node[vert] (n10) at (10,2) {};
 \node[vert] (n11) at (10,4) {};
 \node[vert] (n12) at (3-1.41421356,4+1.41421356) {};

 \draw[-] (n1) to [out=90,in=-132] node[label] {} (n12);
 \draw[-] (n12) to [out=42,in=180] node[label] {} (n2);
 \draw[-] (n2) to [out=0,in=90] node[label] {$[7]$} (n3);
 \draw[-] (n3) to [out=-90,in=0] node[label] {$[6]$} (n4);
 \draw[-] (n4) to [out=180,in=-90] node[label] {$[5]$} (n1);

 \draw[-] (n7) to [out=90,in=180] node[label] {$[3]$} (n8);
 \draw[-] (n8) to [out=0,in=90] node[label] {$[3]$} (n9);
 \draw[-] (n9) to [out=-90,in=0] node[label] {$[4]$} (n10);
 \draw[-] (n10) to [out=180,in=-90] node[label] {$[4]$} (n7);

 \draw[-] (n1) -- node[label][swap] {$[5]$} (n5);
 \draw[-] (n5) -- node[label][swap] {$[6]$} (n3);
 \draw[-] (n7) -- node[label][swap] {$[4]$} (n11);
 \draw[-] (n11) -- node[label] {$[4]$} (n9);

 \draw[-] (n2) to [out=30,in=150] node[label] {$[1]$} (n8);
 \draw[-] (n5) to [out=30,in=150] node[label] {$[2]$} (n11);
 \draw[-] (n4) to [out=-30,in=-150] node[label] {$[2]$} (n10);

\end{tikzpicture}
\caption{The graph $\wG_{[2]}$.}
\end{figure}

Since we have 7 equivalence classes of edges, we may obtain 7 graphs
from $\wG_{[2]}$.  The naming of the new graphs is as follows:
choosing $[1]$ we get $G_7$, choosing $[2]$ we get $G_8$, choosing 
$[3]$ we get $G_5$, choosing $[4]$ we get $G_6$, choosing $[5]$ we get $G_2$,
choosing $[6]$ we get $G_4$, and choosing $[7]$ we get $G_3$.

Thus, there are a total of 14 genus 7 trivalent bridgeless graphs: 
$F_{11},F_{12},F_{13},F_{14}$ and $G_i$ for $i \in \{1, \dots ,10\}$.

\end{longversion} 

\subsection{$H$ has genus 8}
\label{subsec:genus-8}
Since $H$ is cohomology-irreducible, the graphs $H$ and $\wH$ can not be
isomorphic to $E_{42}$: otherwise, $H$ would have genus $\ge10$.

We may choose an edge $e$ so that $H-e$ does not embed into $P$. Since 
$e$ is not a bridge,
we may construct a trivalent graph $\simp(H-e)$ from $H-e$ by
contracting edges which were incident to $e$, as in (3a) or (3b). So
$\simp(H-e)$ is a trivalent graph of genus $7$ which does not embed
into $P$. Hence by our above argument, $\simp(H-e)$ is isomorphic to
one of $F_{11},F_{12},F_{13},F_{14}, G_{i}$ for $i \in \{1, \dots
,10\}$, or a graph $G'$ obtained from $G$ by choosing an edge $e'$,
adding an edge $f$ from the midpoint of $e'$ to an isolated vertex $v$
and then adding a loop $e$ to $v$, as in (c).

In the latter case, $H$ is obtained from the graph $G$ by performing
operation (c) and then (a). But, equivalently, this can be
accomplished by the operations (a) and then (b). 
Thus, to prove
Theorem~\ref{thm:main} for $g=8$, it is sufficient to find $\bZ$-emms
for the finitely many graphs obtained from one of the graphs
$F_{11}$--$F_{14}$, $G_1$--$G_{10}$ by performing one operation of
type (a) or (b).

\section{Genus 6}
\label{sec:6}

In this section, we explain the general method for finding a $\bZ$-emm
for any graph, and illustrate it in the case of the trivalent genus 6
graph $G$.

\subsection{Procedure for a general graph}

Let $\Gamma$ be a directed graph of genus $g$ with edge set $E=\{e_1,
\dots ,e_n\}$. After renaming the edges, we may insist that the edges
$\{e_{g+1}, \dots ,e_n\}$ induce a spanning tree $T$ of~$\Gamma$. Then
for each $e_i$ with $i \in \{1, \dots ,g\}$, we have a corresponding
basis element $f_i$ of the homology group $H_1(\Gamma,\bZ)$, given by:
\begin{displaymath}
  f_i= e_i +\sum_{e_s \in T}b_{i,s}e_s,
  \quad b_{i,s}=0,\pm 1, i\in \{ 1, \dots, g \} ,
\end{displaymath}
and the coedges $e_1^*,\dotsc,e_g^*$ form a basis of the cohomology
group $H^1(G,\bZ)$ (cf. \cite[Lemma 2.3]{AlexeevBrunyate}).

Specifically, $f_i$ is given by the unique simple cycle in
$\Gamma$ which uses only the edge $e_i$ and edges of $T$. 
If we write the vectors $f_i$ as the rows of a $g\times n$ matrix 
then the columns of this matrix are the coedges $e_i^*\in
H^1(G,\bZ)$ written in the basis $\{e_1^*,\dotsc,e_g^*\}$. In
particular, the first $g$ columns form an identity matrix.

Let $q$ be a $\bZ$-emm for $\Gamma$.
Since $q$ is a $\bZ$-valued quadratic form, we may associate to $q$ an
even integral matrix
$M_q=(a_{i,j})$ such that
\begin{displaymath}
  q(x_1, \dots ,x_g)=(x_1, \dots
  ,x_g)\frac{1}{2}M_q(x_1, \dots ,x_g)^T.  
\end{displaymath}
Note here that
$a_{i,j}=a_{j,i}$ is just the coefficient of the term $x_ix_j$ in
$q(x_1, \dots ,x_g)$ if $i \not= j$ and $a_{i,i}$ is just twice the
coefficient of the term $x_i^2$ in $q(x_1, \dots ,x_g)$.

We need to enforce the condition that $q(e^*_i)=1$ for $i=1, \dots
,n$. To ensure that $q(e^*_i)=1$ for $i=1, \dots ,g$ we must have
$a_{i,i}=2$. Now we must ensure that $q(e^*_i)=1$ for $i=g+1, \dots
,n$. This is equivalent to $n-g$ linear equations on $a_{i,j}$:
\begin{displaymath}
  1=\sum_{i=1}^g c_i^2+\sum_{1\leq i<j \leq g}c_ic_ja_{i,j} 
  \quad
  \text{for each column }
  (c_i). 
\end{displaymath}
Further, the condition that $q$ is positive definite implies that each
$a_{i,j}\in \{0,\pm1\}$. Thus, for any given graph, we reduced the
problem to a finite computation.

\subsection{Computation for graph $G$}

We now specialize to graph $G$. In Figure~\ref{fig-g} it is shown as
a labeled directed graph with a
spanning tree denoted by bold edges. 

\begin{showfigures}
\begin{figure}[H]

\caption{The Graph $G$.}
\label{fig-g}
\end{figure}
\end{showfigures}

\noindent Using the spanning tree drawn and the process described
above we get a basis for $H_1(G,\bZ)$, written as the rows of the
following matrix:

\begin{longversion}
$$\begin{array}{l}
  f_1=e_1-e_{12}-e_{13}+e_{14} \\
  f_2=e_2+e_{12}+e_{13}+e_{15} \\
  f_3=e_3+e_8-e_{10}-e_{11}+e_{13}-e_{14} \\
  f_4=e_4+e_7+e_8-e_{11}+e_{13}+e_{15} \\
  f_5=e_5-e_8-e_9+e_{10} \\
  f_6=e_6+e_7+e_8+e_9
  \end{array}$$
\end{longversion}

$$\begin{array}{c|c|c|c|c|c|c|c|c|c|c|c|c|c|c|c}
       & e_1 & e_2 & e_3 & e_4 & e_5 & e_6 & e_7 & e_8 & e_9 & e_{10} & e_{11} & e_{12} & e_{13} & e_{14} & e_{15} \\
\hline
   f_1 & 1 & 0 & 0 & 0 & 0 & 0 & 0 & 0 & 0 & 0 & 0 & -1& -1& 1 & 0\\
   f_2 & 0 & 1 & 0 & 0 & 0 & 0 & 0 & 0 & 0 & 0 & 0 & 1 & 1 & 0 & 1\\
   f_3 & 0 & 0 & 1 & 0 & 0 & 0 & 0 & 1 & 0 & -1& -1& 0 & 1 & -1& 0\\
   f_4 & 0 & 0 & 0 & 1 & 0 & 0 & 1 & 1 & 0 & 0 & -1& 0 & 1 & 0 & 1\\
   f_5 & 0 & 0 & 0 & 0 & 1 & 0 & 0 & -1& -1& 1 & 0 & 0 & 0 & 0 & 0\\
   f_6 & 0 & 0 & 0 & 0 & 0 & 1 & 1 & 1 & 1 & 0 & 0 & 0 & 0 & 0 & 0\\
  \end{array}$$

The linear equations become:
$$\begin{array}{ll}
    (1) & 1=2+a_{4,6} \\
    (2) & 1=4+a_{3,4}-a_{3,5}+a_{3,6}-a_{4,5}+a_{4,6}-a_{5,6} \\
    (3) & 1=2-a_{5,6} \\
    (4) & 1=2-a_{3,5} \\
    (5) & 1=2+a_{3,4} \\
    (6) & 1=2-a_{1,2} \\
    (7) & 1=4-a_{1,2} -a_{1,3}-a_{1,4}+a_{2,3}+a_{2,4}+a_{3,4} \\
    (8) & 1=2-a_{1,3} \\
    (9) & 1=2+a_{2,4} \\
  \end{array}$$
So, equations (1),(3),(4),(5),(6),(8),(9) immediately imply that $1=a_{5,6}=a_{3,5}=a_{1,2}=a_{1,3}$ and $-1=a_{4,6}=a_{3,4}=a_{2,4}$. Applying this
information to (2) and (7) we get $1=a_{3,6}-a_{4,5}$ and $1=a_{2,3}-a_{1,4}$ respectively. Let us arbitrarily choose
$a_{3,6}=a_{2,3}=1$ and $a_{4,5}=a_{1,4}=0$. Hence, we will get a $\bZ$-emm if we can choose the remaining terms of the below matrix in such a way 
that it is positive definite. 
$$\left(\begin{array}{cccccc}
 2 & 1 & 1 & 0 & a_{1,5} & a_{1,6} \\
 1 & 2 & 1 & -1& a_{2,5} & a_{2,6} \\
 1 & 1 & 2 & -1& 1 & 1 \\
 0 & -1& -1& 2 & 0 & -1 \\
 a_{5,1} & a_{5,2} & 1 & 0 & 2 & 1 \\
 a_{6,1} & a_{6,2} & 1 & -1& 1 & 2 \\
  \end{array}\right)$$
One such choice is to set all the unknowns to 0. 
\begin{longversion}
$$\left(\begin{array}{cccccc}
 2 & 1 & 1 & 0 & 0 & 0 \\
 1 & 2 & 1 & -1& 0 & 0 \\
 1 & 1 & 2 & -1& 1 & 1 \\
 0 & -1& -1& 2 & 0 & -1 \\
 0 & 0 & 1 & 0 & 2 & 1 \\
 0 & 0 & 1 & -1& 1 & 2 \\
  \end{array}\right)$$
\end{longversion}
Then the quadratic form corresponding to this matrix is:
$$\begin{array}{rl}
  q(x_1,x_2,x_3,x_4,x_5,x_6)=& x_1^2+x_1x_2+x_1x_3+x_2^2+x_2x_3- \\
			     & x_2x_4+x_3^2-x_3x_4+x_3x_5+x_3x_6+ \\
			     & x_4^2-x_4x_6+x_5^2+x_5x_6+x_6^2
  \end{array}$$

  One can easily check by diagonalizing that this quadratic form is
  indeed positive definite. Moreover, in an appropriately chosen basis,
  it is isomoprhic to the standard quadratic form $E_6$.

\section{Genus 7}
\label{sec:7}

We repeat the general procedure of the previous section for the graphs
$F_{11}$--$F_{14}$ and $G_1$--$G_{10}$. Below, we list one explicit
$\bZ$-emm for each of these graphs. The detailed computations are
available in the long version of this paper.

\begin{showfigures}
\begin{figure}[H]
 %
\caption{The Graph $F_{11}$.}
\label{fig-f11}
\end{figure}
\end{showfigures}

\begin{longversion}
\noindent The basis for $F_{11}$ is given by:
$$\begin{array}{l}
  f_1=e_1+e_{17}-e_{16}-e_{14} \\
  f_2=e_2-e_{15}+e_{16}-e_{17} \\
  f_3=e_3+e_{13}-e_{15}-e_{14} \\
  f_4=e_4-e_{18}-e_{16}+e_{15}-e_{13}+e_{12}+e_{11}-e_{9} \\
  f_5=e_5-e_8-e_{11}-e_{12}+e_{13}-e_{15}+e_{16}+e_{18} \\
  f_6=e_6-e_{10}+e_{11}-e_9 \\
  f_7=e_7-e_8-e_{11}+e_{10} 
  \end{array}$$

Which gives us the following array:
$$\begin{array}{c|c|c|c|c|c|c|c|c|c|c|c|c|c|c|c}
      & e_8 & e_9 & e_{10} & e_{11} & e_{12} & e_{13} & e_{14} & e_{15} & e_{16} & e_{17} & e_{18} \\
\hline
   f_1 & 0 & 0 & 0 & 0 & 0 & 0 & -1 & 0 & -1 & 1 & 0 \\
   f_2 & 0 & 0 & 0 & 0 & 0 & 0 & 0 & -1 & 1 & -1 & 0 \\
   f_3 & 0 & 0 & 0 & 0 & 0 & 1 & -1 & -1 & 0 & 0 & 0\\
   f_4 & 0 & -1 & 0 & 1 & 1 & -1 & 0 & 1 & -1 & 0 & -1\\
   f_5 & -1 & 0 & 0 & -1 & -1 & 1 & 0 & -1 & 1 & 0 & 1 \\
   f_6 & 0 & -1 & -1 & 1 & 0 & 0 & 0 & 0 & 0 & 0 & 0 \\
   f_7 & -1 & 0 & 1 & -1 & 0 & 0 & 0 & 0 & 0 & 0 & 0 \\ 
  \end{array}$$
which gives equations:
$$\begin{array}{ll}
    (1) & 1=2+a_{5,7} \\
    (2) & 1=2+a_{4,6} \\
    (3) & 1=2-a_{6,7} \\
    (4) & 1=4-a_{4,5}+a_{4,6}-a_{4,7}-a_{5,6}+a_{5,7}-a_{6,7} \\
    (5) & 1=2-a_{4,5} \\
    (6) & 1=3-a_{3,4}+a_{3,5}-a_{4,5} \\
    (7) & 1=2+a_{1,3} \\
    (8) & 1=4+a_{2,3}-a_{2,4}+a_{2,5}-a_{3,4}+a_{3,5}-a_{4,5} \\
    (9) & 1=4-a_{1,2}+a_{1,4}-a_{1,5}-a_{2,4}+a_{2,5}-a_{4,5} \\
    (10)& 1=2-a_{1,2} \\
    (11)& 1=2-a_{4,5} \\
  \end{array}$$
The following matrix satisfies these equations and is positive definite:
$$\left(\begin{array}{ccccccc}
    2 & 1& -1& 0 & 0 & -1 & -1 \\
    1 & 2 & 0 & 1& 0& -1 & 0 \\
    -1 & 0 & 2 & 1 & 0 & 0 & 1\\
    0 & 1 & 1 & 2 & 1 & -1 & 0 \\
    0 & 0 & 0 & 1 & 2 & -1 & -1\\
    -1 & -1 & 0 & -1 & -1 & 2 & 1\\
    -1 & 0 & 1 & 0 & -1 & 1 & 2\\
  \end{array}\right)$$
And the quadratic form given by this matrix is:
\end{longversion}
$$%
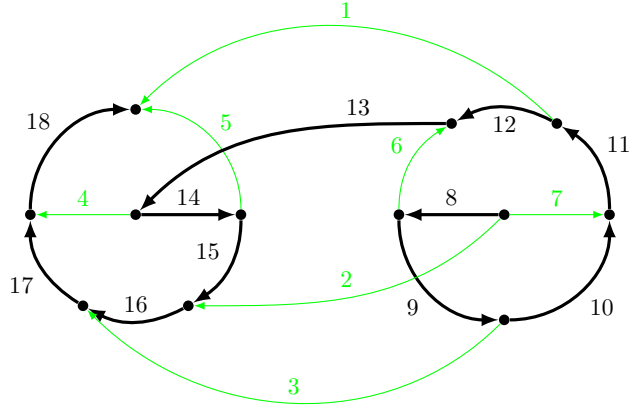
\caption{The graph $F_{12}$.}
\label{fig-f12}
\end{figure}
\end{showfigures}

\begin{longversion}
\noindent The basis for $F_{12}$ is given by:
$$\begin{array}{l}
  f_1=e_1-e_{12}-e_{13}-e_{14}-e_{15}-e_{16}-e_{17}-e_{18} \\
  f_2=e_2-e_8-e_9-e_{10}-e_{11}-e_{12}-e_{13}-e_{14}-e_{15} \\
  f_3=e_3-e_{10}-e_{11}-e_{12}-e_{13}-e_{14}-e_{15}-e_{16} \\
  f_4=e_4-e_{14}-e_{15}-e_{16}-e_{17} \\
  f_5=e_5-e_{15}-e_{16}-e_{17}-e_{18} \\
  f_6=e_6-e_9-e_{10}-e_{11}-e_{12} \\
  f_7=e_7-e_8-e_9-e_{10} 
  \end{array}$$
Which gives us the following array:
$$\begin{array}{c|c|c|c|c|c|c|c|c|c|c|c|c|c|c|c}
      & e_8 & e_9 & e_{10} & e_{11} & e_{12} & e_{13} & e_{14} & e_{15} & e_{16} & e_{17} & e_{18} \\
\hline
   f_1 & 0 & 0 & 0 & 0 & -1 & -1 & -1 & -1 & -1 & -1 & -1 \\
   f_2 & -1 & -1 & -1 & -1 & -1 & -1 & -1 & -1 & 0 & 0 & 0 \\
   f_3 & 0 & 0 & -1 & -1 & -1 & -1 & -1 & -1 & -1 & 0 & 0\\
   f_4 & 0 & 0 & 0 & 0 & 0 & 0 & -1 & -1 & -1 & -1 & 0\\
   f_5 & 0 & 0 & 0 & 0 & 0 & 0 & 0 & -1 & -1 & -1 & -1 \\
   f_6 & 0 & -1 & -1 & -1 & -1 & 0 & 0 & 0 & 0 & 0 & 0 \\
   f_7 & -1 & -1 & -1 & 0 & 0 & 0 & 0 & 0 & 0 & 0 & 0 \\ 
  \end{array}$$
which gives equations:
$$\begin{array}{ll}
    (1) & 1=2+a_{2,7} \\
    (2) & 1=3+a_{2,6}+a_{2,7}+a_{6,7} \\
    (3) & 1=4+a_{2,3}+a_{2,6}+a_{2,7}+a_{3,6}+a_{3,7}+a_{6,7} \\
    (4) & 1=3+a_{2,3}+a_{2,6}+a_{3,6} \\
    (5) & 1=4+a_{1,2}+a_{1,3}+a_{1,6}+a_{2,3}+a_{2,6}+a_{3,6} \\
    (6) & 1=3+a_{1,2}+a_{1,3}+a_{2,3} \\
    (7) & 1=4+a_{1,2}+a_{1,3}+a_{1,4}+a_{2,3}+a_{2,4}+a_{3,4} \\
    (8) & 1=5+a_{1,2}+a_{1,3}+a_{1,4}+a_{1,5}+a_{2,3}+a_{2,4}+a_{2,5}+a_{3,4}+a_{3,5}+a_{4,5} \\
    (9) & 1=4+a_{1,3}+a_{1,4}+a_{1,5}+a_{3,4}+a_{3,5}+a_{4,5} \\
    (10)& 1=3+a_{1,4}+a_{1,5}+a_{4,5} \\
    (11)& 1=2+a_{1,5} \\
  \end{array}$$
The following matrix satisfies these equations and is positive definite:
$$\left(\begin{array}{ccccccc}
    2 & 0& -1& 0 & -1 & 0 & 0 \\
    0 & 2 & -1 & -1& 1& 0 & -1 \\
    -1 & -1 & 2 & 0 & 0 & -1 & 1\\
    0 & -1 & 0 & 2 & -1 & 0 & 0 \\
    -1 & 1 & 0 & -1 & 2 & 0 & 0\\
    0 & 0 & -1 & 0 & 0 & 2 & -1\\
    0 & -1 & 1 & 0 & 0 & -1 & 2\\
  \end{array}\right)$$
And the quadratic form given by this matrix is:
\end{longversion}
$$%
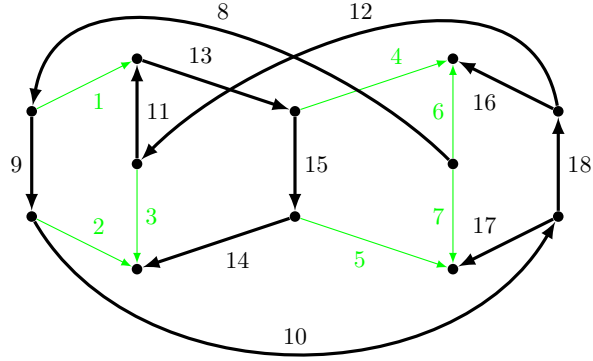
\caption{The Graph $F_{13}$}
\label{fig-f13}
\end{figure}
\end{showfigures}

\begin{longversion}
\noindent The basis for $F_{13}$ is given by:
$$\begin{array}{l}
  f_1=e_1-e_9-e_{10}-e_{11}-e_{12}-e_{18} \\
  f_2=e_2-e_{10}-e_{11}-e_{12}-e_{13}-e_{14}-e_{15}-e_{18} \\
  f_3=e_3-e_{11}-e_{13}-e_{14}-e_{15} \\
  f_4=e_4+e_{11}+e_{12}+e_{13}-e_{16} \\
  f_5=e_5+e_{11}+e_{12}+e_{13}+e_{15}-e_{17}+e_{18} \\
  f_6=e_6-e_8-e_9-e_{10}-e_{16}-e_{18} \\
  f_7=e_7-e_8-e_9-e_{10}-e_{17} 
  \end{array}$$
Which gives us the following array:
$$\begin{array}{c|c|c|c|c|c|c|c|c|c|c|c|c|c|c|c}
      & e_8 & e_9 & e_{10} & e_{11} & e_{12} & e_{13} & e_{14} & e_{15} & e_{16} & e_{17} & e_{18} \\
\hline
   f_1 & 0 & -1 & -1 & -1 & -1 & 0 & 0 & 0 & 0 & 0 & -1 \\
   f_2 & 0 & 0 & -1 & -1 & -1 & -1 & -1 & -1 & 0 & 0 & -1 \\
   f_3 & 0 & 0 & 0 & -1 & 0 & -1 & -1 & -1 & 0 & 0 & 0\\
   f_4 & 0 & 0 & 0 & 1 & 1 & 1 & 0 & 0 & -1 & 0 & 0\\
   f_5 & 0 & 0 & 0 & 1 & 1 & 1 & 0 & 1 & 0 & -1 & 1 \\
   f_6 & -1 & -1 & -1 & 0 & 0 & 0 & 0 & 0 & -1 & 0 & -1 \\
   f_7 & -1 & -1 & -1 & 0 & 0 & 0 & 0 & 0 & 0 & -1 & 0 \\ 
  \end{array}$$
which gives equations:
$$\begin{array}{ll}
    (1) & 1=2+a_{6,7} \\
    (2) & 1=3+a_{1,6}+a_{1,7}+a_{6,7} \\
    (3) & 1=4+a_{1,2}+a_{1,6}+a_{1,7}+a_{2,6}+a_{2,7}+a_{6,7} \\
    (4) & 1=5+a_{1,2}+a_{1,3}-a_{1,4}-a_{1,5}+a_{2,3}-a_{2,4}-a_{2,5}-a_{3,4}-a_{3,5}+a_{4,5} \\
    (5) & 1=4+a_{1,2}-a_{1,4}-a_{1,5}-a_{2,4}-a_{2,5}+a_{4,5} \\
    (6) & 1=4+a_{2,3}-a_{2,4}-a_{2,5}-a_{3,4}-a_{3,5}+a_{4,5} \\
    (7) & 1=2+a_{2,3} \\
    (8) & 1=3+a_{2,3}-a_{2,5}-a_{3,5} \\
    (9) & 1=2+a_{4,6} \\
    (10)& 1=2+a_{5,7} \\
    (11)& 1=4+a_{1,2}-a_{1,5}+a_{1,6}-a_{2,5}+a_{2,6}-a_{5,6} \\
  \end{array}$$
The following matrix satisfies these equations and is positive definite:
$$\left(\begin{array}{ccccccc}
    2 & 0& 0& 0 & 1 & 0 & -1 \\
    0 & 2 & -1 & 0& 1& -1 & 0 \\
    0 & -1 & 2 & 0 & 0 & 0 & 0\\
    0 & 0 & 0 & 2 & -1 & -1 & 1 \\
    1 & 1 & 0 & -1 & 2 & 0 & -1\\
    0 & -1 & 0 & -1 & 0 & 2 & -1\\
    -1 & 0 & 0 & 1 & -1 & -1 & 2\\
  \end{array}\right)$$
And the quadratic form given by this matrix is:
\end{longversion}
$$%
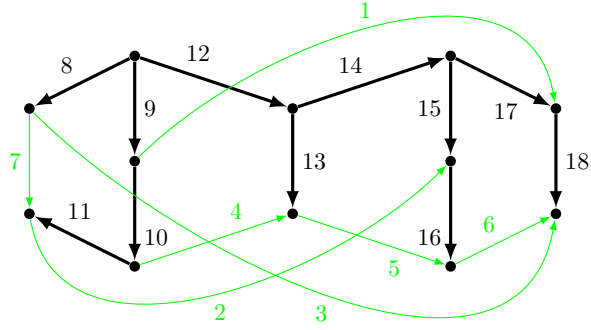
\caption{The Graph $F_{14}$}
\label{fig-f14}
\end{figure}
\end{showfigures}

\begin{longversion}
\noindent The basis for $F_{14}$ is given by:
$$\begin{array}{l}
  f_1=e_1+e_9-e_{12}-e_{14}-e_{17} \\
  f_2=e_2+e_9+e_{10}+e_{11}-e_{12}-e_{14}-e_{15} \\
  f_3=e_3+e_8-e_{12}-e_{14}-e_{17}-e_{18} \\
  f_4=e_4+e_9+e_{10}-e_{12}-e_{13} \\
  f_5=e_5+e_{13}-e_{14}-e_{15}-e_{16} \\
  f_6=e_6+e_{15}+e_{16}-e_{17}-e_{18} \\
  f_7=e_7+e_8-e_9-e_{10}-e_{11} \\
  \end{array}$$
Which gives us the following array:
$$\begin{array}{c|c|c|c|c|c|c|c|c|c|c|c|c|c|c|c}
       & e_8 & e_9 & e_{10} & e_{11} & e_{12} & e_{13} & e_{14} & e_{15} & e_{16} & e_{17} & e_{18} \\
\hline
   f_1 & 0 & 1 & 0 & 0 & -1 & 0 & -1 & 0 & 0 & -1 & 0 \\
   f_2 & 0 & 1 & 1 & 1 & -1 & 0 & -1 & -1 & 0 & 0 & 0 \\
   f_3 & 1 & 0 & 0 & 0 & -1 & 0 & -1 & 0 & 0 & -1 & -1\\
   f_4 & 0 & 1 & 1 & 0 & -1 & -1 & 0 & 0 & 0 & 0 & 0\\
   f_5 & 0 & 0 & 0 & 0 & 0 & 1 & -1 & -1& -1& 0 & 0 \\
   f_6 & 0 & 0 & 0 & 0 & 0 & 0 & 0 & 1 & 1 & -1 & -1 \\
   f_7 & 1 & -1 & -1 & -1 & 0 & 0 & 0 & 0 & 0 & 0 & 0 \\ 
  \end{array}$$
which gives equations:
$$\begin{array}{ll}
    (1) & 1=2+a_{3,7} \\
    (2) & 1=4+a_{1,2}+a_{1,4}-a_{1,7}+a_{2,4}-a_{2,7}-a_{4,7} \\
    (3) & 1=3+a_{2,4}-a_{2,7}-a_{4,7} \\
    (4) & 1=2-a_{2,7} \\
    (5) & 1=4+a_{1,2}+a_{1,3}+a_{1,4}+a_{2,3}+a_{2,4}+a_{3,4} \\
    (6) & 1=2-a_{4,5} \\
    (7) & 1=4+a_{1,2}+a_{1,3}+a_{1,5}+a_{2,3}+a_{2,5}+a_{3,5} \\
    (8) & 1=3+a_{2,5}-a_{2,6}-a_{5,6} \\
    (9) & 1=2-a_{5,6} \\
    (10)& 1=3+a_{1,3}+a_{1,6}+a_{3,6} \\
    (11)& 1=2+a_{3,6} \\
  \end{array}$$
The following matrix satisfies these equations and is positive definite:
$$\left(\begin{array}{ccccccc}
    2 & -1& -1& 0 & 0 & 0 & 0 \\
    -1& 2 & 0 & -1& -1& 0 & 1 \\
    -1& 0 & 2 & 0 & 0 & -1& -1\\
    0 & -1& 0 & 2 & 1 & 0 & 0 \\
    0 & -1& 0 & 1 & 2 & 1 & -1\\
    0 & 0 & -1 & 0 & 1 & 2 & 0\\
    0 & 1 & -1 & 0 & -1 & 0 & 2\\
  \end{array}\right)$$
And the quadratic form given by this matrix is:
\end{longversion}
$$%
\caption{The graph $G_1$.}
\label{fig-g1}
\end{figure}
\end{showfigures}

\begin{longversion}
\noindent The basis for $G_1$ is given by:
$$\begin{array}{l}
  f_1=e_1+e_{12}+e_{13}+e_{16}+e_{18} \\
  f_2=e_2+e_{12}+e_{15}+e_{16}+e_{18} \\
  f_3=e_3-e_{13}-e_{14}+e_{15} \\
  f_4=e_4+e_8+e_9+e_{11}+e_{12}+e_{13}+e_{14} \\
  f_5=e_5-e_8-e_9+e_{10} \\
  f_6=e_6-e_{10}-e_{11}-e_{17}+e_{18} \\
  f_7=e_7+e_{9}+e_{11}+e_{17}-e_{18}  
  \end{array}$$
Which gives us the following array:
$$\begin{array}{c|c|c|c|c|c|c|c|c|c|c|c|c|c|c|c}
       & e_8 & e_9 & e_{10} & e_{11} & e_{12} & e_{13} & e_{14} & e_{15} & e_{16} & e_{17} & e_{18} \\
\hline
   f_1 & 0 & 0 & 0 & 0 & 1 & 1 & 0 & 0 & 1 & 0 & 1 \\
   f_2 & 0 & 0 & 0 & 0 & 1 & 0 & 0 & 1 & 1 & 0 & 1 \\
   f_3 & 0 & 0 & 0 & 0 & 0 & -1 & -1 & 1 & 0 & 0 & 0\\
   f_4 & 1 & 1 & 0 & 1 & 1 & 1 & 1 & 0 & 0 & 0 & 0\\
   f_5 & -1 & -1 & 1 & 0 & 0 & 0 & 0 & 0 & 0 & 0 & 0 \\
   f_6 & 0 & 0 & -1 & -1 & 0 & 0 & 0 & 0 & 0 & -1 & 1 \\
   f_7 & 0 & 1 & 0 & 1 & 0 & 0 & 0 & 0 & 0 & 1 & -1 \\ 
  \end{array}$$
which gives equations:
$$\begin{array}{ll}
    (1) & 1=2-a_{4,5} \\
    (2) & 1=3-a_{4,5}+a_{4,7}-a_{5,7} \\
    (3) & 1=2-a_{5,6} \\
    (4) & 1=3-a_{4,6}+a_{4,7}-a_{6,7} \\
    (5) & 1=3+a_{1,2}+a_{1,4}+a_{2,4} \\
    (6) & 1=3-a_{1,3}+a_{1,4}-a_{3,4} \\
    (7) & 1=2-a_{3,4} \\
    (8) & 1=2+a_{2,3} \\
    (9) & 1=2+a_{1,2} \\
    (10)& 1=2-a_{6,7} \\
    (11)& 1=4+a_{1,2}+a_{1,6}-a_{1,7}+a_{2,6}-a_{2,7}-a_{6,7} \\
  \end{array}$$
The following matrix satisfies these equations and is positive definite:
$$\left(\begin{array}{ccccccc}
    2 & -1 & 1 & 0 & 1 & 1 & 1 \\
    -1 & 2 & -1 & -1 & -1& -1 & 0 \\
    1& -1 & 2 & 1 & 1 & 1& 0\\
    0 & -1& 1 & 2 & 1 & 0 & -1 \\
    1 & -1& 1 & 1 & 2 & 1 & 0\\
    1 & -1 & 1 & 0 & 1 & 2 & 1\\
    1 & 0 & 0 & -1 & 0 & 1 & 2\\
  \end{array}\right)$$
And the quadratic form given by this matrix is:
\end{longversion}
$$%
\caption{The graph $G_2$.}
\label{fig-g2}
\end{figure}
\end{showfigures}

\begin{longversion}
\noindent The basis for $G_2$ is given by:
$$\begin{array}{l}
  f_1=e_1-e_8+e_9 \\
  f_2=e_2-e_9-e_{10}-e_{11}+e_{13} \\
  f_3=e_3-e_8-e_{10}-e_{11}+e_{14} \\
  f_4=e_4+e_{11}-e_{12}-e_{14}-e_{15}+e_{17} \\
  f_5=e_5+e_{15}-e_{17}-e_{18} \\
  f_6=e_6+e_{16}-e_{17}-e_{18} \\
  f_7=e_7+e_{11}-e_{12}-e_{13}-e_{16}+e_{17} 
  \end{array}$$
Which gives us the following array:
$$\begin{array}{c|c|c|c|c|c|c|c|c|c|c|c|c|c|c|c}
       & e_8 & e_9 & e_{10} & e_{11} & e_{12} & e_{13} & e_{14} & e_{15} & e_{16} & e_{17} & e_{18} \\
\hline
   f_1 & -1 & 1 & 0 & 0 & 0 & 0 & 0 & 0 & 0 & 0 & 0 \\
   f_2 & 0 & -1 & -1 & -1 & 0 & 1 & 0 & 0 & 0 & 0 & 0 \\
   f_3 & -1 & 0 & -1 & -1 & 0 & 0 & 1 & 0 & 0 & 0 & 0\\
   f_4 & 0 & 0 & 0 & 1 & -1 & 0 & -1 & -1 & 0 & 1 & 0\\
   f_5 & 0 & 0 & 0 & 0 & 0 & 0 & 0 & 1 & 0 & -1 & -1 \\
   f_6 & 0 & 0 & 0 & 0 & 0 & 0 & 0 & 0 & 1 & -1 & -1 \\
   f_7 & 0 & 0 & 0 & 1 & -1 & -1 & 0 & 0 & -1 & 1 & 0 \\ 
  \end{array}$$
which gives equations:
$$\begin{array}{ll}
    (1) & 1=2+a_{1,3} \\
    (2) & 1=2-a_{1,2} \\
    (3) & 1=2+a_{2,3} \\
    (4) & 1=4+a_{2,3}-a_{2,4}-a_{2,7}-a_{3,4}-a_{3,7}+a_{4,7} \\
    (5) & 1=2+a_{4,7} \\
    (6) & 1=2-a_{2,7} \\
    (7) & 1=2-a_{3,4} \\
    (8) & 1=2-a_{4,5} \\
    (9) & 1=2-a_{6,7} \\
    (10)& 1=4-a_{4,5}-a_{4,6}+a_{4,7}+a_{5,6}-a_{5,7}-a_{6,7} \\
    (11)& 1=2+a_{5,6} \\
  \end{array}$$
The following matrix satisfies these equations and is positive definite:
$$\left(\begin{array}{ccccccc}
    2 & 1 & -1& 0 & 0 & 0 & 0 \\
    1 & 2 & -1 & -1 & 0& 0 & 1 \\
    -1& -1 & 2 & 1 & 0 & 0& 0\\
    0 & -1& 1 & 2 & 1 & -1 & -1 \\
    0 & 0& 0 & 1 & 2 & -1 & 0\\
    0 & 0 & 0 & -1 & -1 & 2 & 1\\
    0 & 1 & 0 & -1 & 0 & 1 & 2\\
  \end{array}\right)$$
And the quadratic form given by this matrix is:
\end{longversion}
$$%
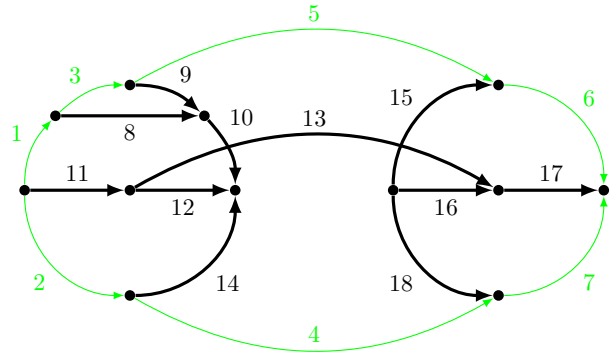
\caption{The Graph $G_3$}
\label{fig-g3}
\end{figure}
\end{showfigures}

\begin{longversion}
\noindent The basis for $G_3$ is given by:
$$\begin{array}{l}
  f_1=e_1+e_8+e_{10}-e_{11}-e_{12} \\
  f_2=e_2-e_{11}-e_{12}+e_{14} \\
  f_3=e_3-e_8+e_9 \\
  f_4=e_4+e_{12}-e_{13}-e_{14}+e_{16}-e_{18} \\
  f_5=e_5-e_9-e_{10}+e_{12}-e_{13}-e_{15}+e_{16} \\
  f_6=e_6+e_{15}-e_{16}-e_{17} \\
  f_7=e_7-e_{16}-e_{17}+e_{18} \\
  \end{array}$$
Which gives us the following array:
$$\begin{array}{c|c|c|c|c|c|c|c|c|c|c|c|c|c|c|c}
       & e_8 & e_9 & e_{10} & e_{11} & e_{12} & e_{13} & e_{14} & e_{15} & e_{16} & e_{17} & e_{18} \\
\hline
   f_1 & 1 & 0 & 1 & -1 & -1 & 0 & 0 & 0 & 0 & 0 & 0 \\
   f_2 & 0 & 0 & 0 & -1 & -1 & 0 & 1 & 0 & 0 & 0 & 0 \\
   f_3 & -1 & 1 & 0 & 0 & 0 & 0 & 0 & 0 & 0 & 0 & 0\\
   f_4 & 0 & 0 & 0 & 0 & 1 & -1 & -1 & 0 & 1 & 0 & -1\\
   f_5 & 0 & -1 & -1 & 0 & 1 & -1 & 0 & -1 & 1 & 0 & 0 \\
   f_6 & 0 & 0 & 0 & 0 & 0 & 0 & 0 & 1 & -1 & -1 & 0 \\
   f_7 & 0 & 0 & 0 & 0 & 0 & 0 & 0 & 0 & -1 & -1 & 1 \\ 
  \end{array}$$
which gives equations:
$$\begin{array}{ll}
    (1) & 1=2-a_{1,3} \\
    (2) & 1=2-a_{3,5} \\
    (3) & 1=2-a_{1,5} \\
    (4) & 1=2+a_{1,2} \\
    (5) & 1=4+a_{1,2}-a_{1,4}-a_{1,5}-a_{2,4}-a_{2,5}+a_{4,5} \\
    (6) & 1=2+a_{4,5} \\
    (7) & 1=2-a_{2,4} \\
    (8) & 1=2-a_{5,6} \\
    (9) & 1=4+a_{4,5}-a_{4,6}-a_{4,7}-a_{5,6}-a_{5,7}+a_{6,7} \\
    (10)& 1=2+a_{6,7} \\
    (11)& 1=2-a_{4,7} \\
  \end{array}$$
The following matrix satisfies these equations and is positive definite:
$$\left(\begin{array}{ccccccc}
    2 & -1 & 1 & -1 & 1 & 0 & 0 \\
    -1 & 2 & 0 & 1 & 0 & 0 & 0 \\
    1& 0 & 2 & 0 & 1 & 0& 0\\
    -1 & 1 & 0 & 2 & -1 & -1 & 1 \\
    1 & 0& 1 & -1 & 2 & 1 & 0\\
    0 & 0 & 0 & -1 & 1 & 2 & -1\\
    0 & 0 & 0 & 1 & 0 & -1 & 2\\
  \end{array}\right)$$
And the quadratic form given by this matrix is:
\end{longversion}
$$%
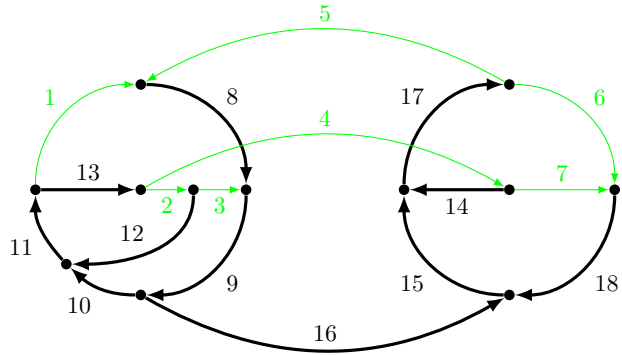
\caption{The Graph $G_4$}
\label{fig-g4}
\end{figure}
\end{showfigures}

\begin{longversion}
\noindent The basis for $G_4$ is given by:
$$\begin{array}{l}
  f_1=e_1+e_8+e_9+e_{10}+e_{11} \\
  f_2=e_2+e_{11}+e_{12}+e_{13} \\
  f_3=e_3+e_9+e_{10}-e_{12} \\
  f_4=e_4+e_{10}+e_{11}+e_{13}+e_{14}-e_{15}-e_{16} \\
  f_5=e_5+e_8+e_9+e_{15}+e_{16}+e_{17} \\
  f_6=e_6+e_{15}+e_{17}+e_{18} \\
  f_7=e_7-e_{14}+e_{15}+e_{18} 
  \end{array}$$
Which gives us the following array:
$$\begin{array}{c|c|c|c|c|c|c|c|c|c|c|c|c|c|c|c}
        & e_8 & e_9 & e_{10} & e_{11} & e_{12} & e_{13} & e_{14} & e_{15} & e_{16} & e_{17} & e_{18} \\
\hline
   f_1 & 1 & 1 & 1 & 1 & 0 & 0 & 0 & 0 & 0 & 0 & 0 \\
   f_2 & 0 & 0 & 0 & 1 & 1 & 1 & 0 & 0 & 0 & 0 & 0 \\
   f_3 & 0 & 1 & 1 & 0 & -1 & 0 & 0 & 0 & 0 & 0 & 0 \\
   f_4 & 0 & 0 & 1 & 1 & 0 & 1 & 1 & -1 & -1 & 0 & 0 \\
   f_5 & 1 & 1 & 0 & 0 & 0 & 0 & 0 & 1 & 1 & 1 & 0 \\
   f_6 & 0 & 0 & 0 & 0 & 0 & 0 & 0 & 1 & 0 & 1 & 1 \\
   f_7 & 0 & 0 & 0 & 0 & 0 & 0 & -1 & 1 & 0 & 0 & 1 \\ 
  \end{array}$$
which gives equations:
$$\begin{array}{ll}
    (1) & 1=2+a_{1,5} \\
    (2) & 1=3+a_{1,3}+a_{1,5}+a_{3,5} \\
    (3) & 1=3+a_{1,3}+a_{1,4}+a_{3,4} \\
    (4) & 1=3+a_{1,2}+a_{1,4}+a_{2,4} \\
    (5) & 1=2-a_{2,3} \\
    (6) & 1=2+a_{2,4} \\
    (7) & 1=2-a_{4,7} \\
    (8) & 1=4-a_{4,5}-a_{4,6}-a_{4,7}+a_{5,6}+a_{5,7}+a_{6,7} \\
    (9) & 1=2-a_{4,5} \\
    (10)& 1=2+a_{5,6} \\
    (11)& 1=2+a_{6,7} \\
  \end{array}$$
The following matrix satisfies these equations and is positive definite:
$$\left(\begin{array}{ccccccc}
    2 & 0 & 0 & -1 & -1 & 0 & -1 \\
    0 & 2 & 1 & -1 & -1 & 0 & 0 \\
    0 & 1 & 2 & -1 & -1 & 0 & 0 \\
    -1 & -1 & -1 & 2 & 1 & 0 & 1 \\
    -1 & -1 & -1 & 1 & 2 & -1 & 1\\
    0 & 0 & 0 & 0 & -1 & 2 & -1\\
    -1 & 0 & 0 & 1 & 1 & -1 & 2\\
  \end{array}\right)$$
And the quadratic form given by this matrix is:
\end{longversion}
$$%
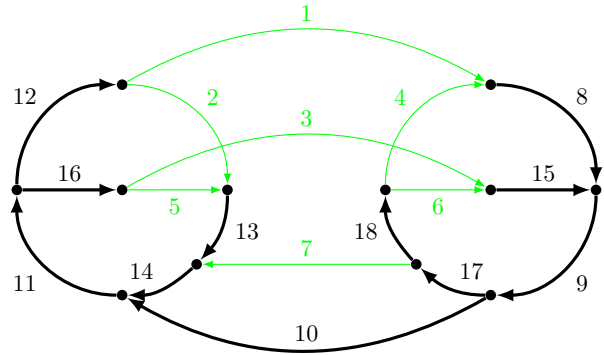
\caption{The Graph $G_5$}
\label{fig-g5}
\end{figure}
\end{showfigures}

\begin{longversion}
\noindent The basis for $G_5$ is given by:
$$\begin{array}{l}
  f_1=e_1+e_8+e_9+e_{10}+e_{11}+e_{12} \\
  f_2=e_2+e_{11}+e_{12}+e_{13}+e_{14} \\
  f_3=e_3+e_9+e_{10}+e_{11}+e_{15}+e_{16} \\
  f_4=e_4+e_8+e_9+e_{17}+e_{18} \\
  f_5=e_5+e_{11}+e_{13}+e_{14}+e_{16} \\
  f_6=e_6+e_9+e_{15}+e_{17}+e_{18} \\
  f_7=e_7-e_{10}+e_{14}+e_{17} \\
  \end{array}$$
Which gives us the following array:
$$\begin{array}{c|c|c|c|c|c|c|c|c|c|c|c|c|c|c|c}
        & e_8 & e_9 & e_{10} & e_{11} & e_{12} & e_{13} & e_{14} & e_{15} & e_{16} & e_{17} & e_{18} \\
\hline
   f_1 & 1 & 1 & 1 & 1 & 1 & 0 & 0 & 0 & 0 & 0 & 0 \\
   f_2 & 0 & 0 & 0 & 1 & 1 & 1 & 1 & 0 & 0 & 0 & 0 \\
   f_3 & 0 & 1 & 1 & 1 & 0 & 0 & 0 & 1 & 1 & 0 & 0 \\
   f_4 & 1 & 1 & 0 & 0 & 0 & 0 & 0 & 0 & 0 & 1 & 1 \\
   f_5 & 0 & 0 & 0 & 1 & 0 & 1 & 1 & 0 & 1 & 0 & 0 \\
   f_6 & 0 & 1 & 0 & 0 & 0 & 0 & 0 & 1 & 0 & 1 & 1 \\
   f_7 & 0 & 0 & -1& 0 & 0 & 0 & 1 & 0 & 0 & 1 & 0 \\ 
  \end{array}$$
which gives equations:
$$\begin{array}{ll}
    (1) & 1=2+a_{1,4} \\
    (2) & 1=4+a_{1,3}+a_{1,4}+a_{1,6}+a_{3,4}+a_{3,6}+a_{4,6} \\
    (3) & 1=3+a_{1,3}-a_{1,7}-a_{3,7} \\
    (4) & 1=4+a_{1,2}+a_{1,3}+a_{1,5}+a_{2,3}+a_{2,5}+a_{3,5} \\
    (5) & 1=2+a_{1,2} \\
    (6) & 1=2+a_{2,5} \\
    (7) & 1=3+a_{2,5}+a_{2,7}+a_{5,7} \\
    (8) & 1=2+a_{3,6} \\
    (9) & 1=2+a_{3,5} \\
    (10)& 1=3+a_{4,6}+a_{4,7}+a_{6,7} \\
    (11)& 1=2+a_{4,6} \\
  \end{array}$$
The following matrix satisfies these equations and is positive definite:
$$\left(\begin{array}{ccccccc}
    2 & -1& 0 & -1 & 0 & 0 & 1 \\
    -1& 2 & 0 & 0 & -1 & 0 & 0 \\
    0 & 0 & 2 & 0 & -1 & -1 & 1 \\
    -1 & 0 & 0 & 2 & 0 & -1 & 0 \\
    0 & -1 & -1 & 0 & 2 & 1 & -1\\
    0 & 0 & -1 & -1 & 1 & 2 & -1\\
    1 & 0 & 1 & 0 & -1 & -1 & 2\\
  \end{array}\right)$$
And the quadratic form given by this matrix is:
\end{longversion}
$$%
\caption{The Graph $G_6$}
\label{fig-g6}
\end{figure}
\end{showfigures}

\begin{longversion}
\noindent The basis for $G_6$ is given by:
$$\begin{array}{l}
  f_1=e_1-e_8-e_9-e_{12}-e_{14}-e_{18} \\
  f_2=e_2-e_8-e_{11}+e_{13}+e_{15}+e_{17} \\
  f_3=e_3-e_8+e_{10}-e_{12}+e_{16}-e_{18} \\
  f_4=e_4+e_9+e_{12}+e_{13} \\
  f_5=e_5+e_{10}+e_{11}-e_{12}-e_{13} \\
  f_6=e_6+e_{14}+e_{15}+e_{17}+e_{18} \\
  f_7=e_7+e_{16}-e_{17}-e_{18} \\
  \end{array}$$
Which gives us the following array:
$$\begin{array}{c|c|c|c|c|c|c|c|c|c|c|c|c|c|c|c}
        & e_8 & e_9 & e_{10} & e_{11} & e_{12} & e_{13} & e_{14} & e_{15} & e_{16} & e_{17} & e_{18} \\
\hline
   f_1 & -1 & -1 & 0 & 0 & -1 & 0 & -1 & 0 & 0 & 0 & -1 \\
   f_2 & -1 & 0 & 0 & -1 & 0 & 1 & 0 & 1 & 0 & 1 & 0 \\
   f_3 & -1 & 0 & 1 & 0 & -1 & 0 & 0 & 0 & 1 & 0 & -1 \\
   f_4 & 0 & 1 & 0 & 0 & 1 & 1 & 0 & 0 & 0 & 0 & 0 \\
   f_5 & 0 & 0 & 1 & 1 & -1 & -1 & 0 & 0 & 0 & 0 & 0 \\
   f_6 & 0 & 0 & 0 & 0 & 0 & 0 & 1 & 1 & 0 & 1 & 1 \\
   f_7 & 0 & 0 & 0 & 0 & 0 & 0 & 0 & 0 & 1 & -1 & -1 \\ 
  \end{array}$$
which gives equations:
$$\begin{array}{ll}
    (1) & 1=3+a_{1,2}+a_{1,3}+a_{2,3} \\
    (2) & 1=2-a_{1,4} \\
    (3) & 1=2+a_{3,5} \\
    (4) & 1=2-a_{2,5} \\
    (5) & 1=4+a_{1,3}-a_{1,4}+a_{1,5}-a_{3,4}+a_{3,5}-a_{4,5} \\
    (6) & 1=3+a_{2,4}-a_{2,5}-a_{4,5} \\
    (7) & 1=2-a_{1,6} \\
    (8) & 1=2+a_{2,6} \\
    (9) & 1=2+a_{3,7} \\
    (10)& 1=3+a_{2,6}-a_{2,7}-a_{6,7} \\
    (11)& 1=4+a_{1,3}-a_{1,6}+a_{1,7}-a_{3,6}+a_{3,7}-a_{6,7} \\
  \end{array}$$
The following matrix satisfies these equations and is positive definite:
$$\left(\begin{array}{ccccccc}
    2 & -1 & -1 & 1 & 0 & 1 & 0 \\
    -1 & 2 & 0 & -1 & 1 & -1 & 0 \\
    -1 & 0 & 2 & 0 & -1 & -1 & -1 \\
    1 & -1 & 0 & 2 & 0 & 0 & -1 \\
    0 & 1 & -1 & 0 & 2 & 0 & 0\\
    1 & -1 & -1 & 0 & 0 & 2 & 1\\
    0 & 0 & -1 & -1 & 0 & 1 & 2\\
  \end{array}\right)$$
And the quadratic form given by this matrix is:
\end{longversion}
$$%
\caption{The Graph $G_7$}
\label{fig-g7}
\end{figure}
\end{showfigures}

\begin{longversion}
\noindent The basis for $G_7$ is given by:
$$\begin{array}{l}
  f_1=e_1-e_8+e_{12}-e_{15}-e_{17}+e_{18} \\
  f_2=e_2-e_8+e_9+e_{12}-e_{13}-e_{15}+e_{16} \\
  f_3=e_3-e_8+e_9-e_{10}+e_{11}-e_{15}+e_{16} \\
  f_4=e_4+e_{10}+e_{12}+e_{18} \\
  f_5=e_5-e_{10}+e_{11}-e_{12}+e_{13} \\
  f_6=e_6+e_{14}+e_{15}+e_{17} \\
  f_7=e_7-e_{14}-e_{15}+e_{16} \\
  \end{array}$$
Which gives us the following array:
$$\begin{array}{c|c|c|c|c|c|c|c|c|c|c|c|c|c|c|c}
        & e_8 & e_9 & e_{10} & e_{11} & e_{12} & e_{13} & e_{14} & e_{15} & e_{16} & e_{17} & e_{18} \\
\hline
   f_1 & -1 & 0 & 0  & 0 & 1  & 0  & 0  & -1 & 0 & -1 & 1 \\
   f_2 & -1 & 1 & 0  & 0 & 1  & -1 & 0  & -1 & 1 & 0  & 0 \\
   f_3 & -1 & 1 & -1 & 1 & 0  & 0  & 0  & -1 & 1 & 0  & 0 \\
   f_4 & 0  & 0 & 1  & 0 & 1  & 0  & 0  & 0  & 0 & 0  & 1 \\
   f_5 & 0  & 0 & -1 & 1 & -1 & 1  & 0  & 0  & 0 & 0  & 0 \\
   f_6 & 0  & 0 & 0  & 0 & 0  & 0  & 1  & 1  & 0 & 1  & 0 \\
   f_7 & 0  & 0 & 0  & 0 & 0  & 0  & -1 & -1 & 1 & 0  & 0 \\ 
  \end{array}$$
which gives equations:
$$\begin{array}{ll}
    (1) & 1=3+a_{1,2}+a_{1,3}+a_{2,3} \\
    (2) & 1=2+a_{2,3} \\
    (3) & 1=3-a_{3,4}+a_{3,5}-a_{4,5} \\
    (4) & 1=2+a_{3,5} \\
    (5) & 1=4+a_{1,2}+a_{1,4}-a_{1,5}+a_{2,4}-a_{2,5}-a_{4,5} \\
    (6) & 1=2-a_{2,5} \\
    (7) & 1=2-a_{6,7} \\
    (8) & 1=5+a_{1,2}+a_{1,3}-a_{1,6}+a_{1,7}+a_{2,3}-a_{2,6}+a_{2,7}-a_{3,6}+a_{3,7}-a_{6,7} \\
    (9) & 1=3+a_{2,3}+a_{2,7}+a_{3,7} \\
    (10)& 1=2-a_{1,6} \\
    (11)& 1=2+a_{1,4} \\
  \end{array}$$
The following matrix satisfies these equations and is positive definite:
$$\left(\begin{array}{ccccccc}
    2 & -1 & 0 & -1 & -1 & 1 & 0 \\
    -1 & 2 & -1 & 0 & 1 & -1 & -1 \\
    0 & -1 & 2 & 0 & -1 & 0 & 0 \\
    -1 & 0 & 0 & 2 & 1 & -1 & 0 \\
    -1 & 1 & -1 & 1 & 2 & -1 & 0\\
    1 & -1 & 0 & -1 & -1 & 2 & 1\\
    0 & -1 & 0 & 0 & 0 & 1 & 2\\
  \end{array}\right)$$
And the quadratic form given by this matrix is:
\end{longversion}
$$%
\caption{The Graph $G_8$}
\label{fig-g8}
\end{figure}
\end{showfigures}

\begin{longversion}
\noindent The basis for $G_8$ is given by:
$$\begin{array}{l}
  f_1=e_1-e_8-e_{10}-e_{13}+e_{14} \\
  f_2=e_2-e_8-e_{10}+e_{12}-e_{13} \\
  f_3=e_3-e_9+e_{10}+e_{13}-e_{14} \\
  f_4=e_4+e_9-e_{10}-e_{11}-e_{15}+e_{17} \\
  f_5=e_5-e_{11}-e_{12}+e_{13}-e_{16}+e_{17} \\
  f_6=e_6+e_{16}-e_{17}-e_{18} \\
  f_7=e_7+e_{15}-e_{17}-e_{18} \\
  \end{array}$$
Which gives us the following array:
$$\begin{array}{c|c|c|c|c|c|c|c|c|c|c|c|c|c|c|c}
        & e_8 & e_9 & e_{10} & e_{11} & e_{12} & e_{13} & e_{14} & e_{15} & e_{16} & e_{17} & e_{18} \\
\hline
   f_1 & -1& 0 & -1& 0 & 0 & -1& 1 & 0 & 0 & 0 & 0 \\
   f_2 & -1& 0 & -1& 0 & 1 & -1& 0 & 0 & 0 & 0 & 0 \\
   f_3 & 0 & -1& 1 & 0 & 0 & 1 & -1& 0 & 0 & 0 & 0 \\
   f_4 & 0 & 1 & -1& -1& 0 & 0 & 0 & -1& 0 & 1 & 0 \\
   f_5 & 0 & 0 & 0 & -1& -1& 1 & 0 & 0 & -1& 1 & 0 \\
   f_6 & 0 & 0 & 0 & 0 & 0 & 0 & 0 & 0 & 1 & -1& -1 \\
   f_7 & 0 & 0 & 0 & 0 & 0 & 0 & 0 & 1 & 0 & -1& -1 \\ 
  \end{array}$$
which gives equations:
$$\begin{array}{ll}
    (1) & 1=2+a_{1,2} \\
    (2) & 1=2-a_{3,4} \\
    (3) & 1=4+a_{1,2}-a_{1,3}+a_{1,4}-a_{2,3}+a_{2,4}-a_{3,4} \\
    (4) & 1=2+a_{4,5} \\
    (5) & 1=2-a_{2,5} \\
    (6) & 1=4+a_{1,2}-a_{1,3}-a_{1,5}-a_{2,3}-a_{2,5}+a_{3,5} \\
    (7) & 1=2-a_{1,3} \\
    (8) & 1=2-a_{4,7} \\
    (9) & 1=2-a_{5,6} \\
    (10)& 1=4+a_{4,5}-a_{4,6}-a_{4,7}-a_{5,6}-a_{5,7}+a_{6,7} \\
    (11)& 1=2+a_{6,7} \\
  \end{array}$$
The following matrix satisfies these equations and is positive definite:
$$\left(\begin{array}{ccccccc}
    2 & -1& 1 & 0 & 0 & 0 & 0 \\
    -1& 2 & 0 & 0 & 1 & 0 & 0 \\
    1 & 0 & 2 & 1 & 0 & 0 & 0 \\
    0 & 0 & 1 & 2 & -1 & -1 & 1 \\
    0 & 1 & 0 & -1 & 2 & 1 & 0\\
    0 & 0 & 0 & -1 & 1 & 2 & -1\\
    0 & 0 & 0 & 1 & 0 & -1 & 2\\
  \end{array}\right)$$
And the quadratic form given by this matrix is:
\end{longversion}
$$%
\caption{The Graph $G_9$}
\label{fig-g9}
\end{figure}
\end{showfigures}

\begin{longversion}
\noindent The basis for $G_9$ is given by:
$$\begin{array}{l}
  f_1=e_1+e_8+e_9+e_{10} \\
  f_2=e_2-e_8-e_9+e_{11} \\
  f_3=e_3-e_{13} \\
  f_4=e_4-e_9+e_{11}-e_{12}-e_{13}-e_{14}-e_{15}+e_{16} \\
  f_5=e_5+e_{10}+e_{11}-e_{12}-e_{13}-e_{14}-e_{15}-e_{18} \\
  f_6=e_6+e_{16}+e_{17}+e_{18} \\
  f_7=e_7+e_{15}-e_{16}-e_{17} 
  \end{array}$$
Which gives us the following array:
$$\begin{array}{c|c|c|c|c|c|c|c|c|c|c|c|c|c|c|c}
        & e_8 & e_9 & e_{10} & e_{11} & e_{12} & e_{13} & e_{14} & e_{15} & e_{16} & e_{17} & e_{18} \\
\hline
   f_1 & 1 & 1 & 1 & 0 & 0 & 0 & 0 & 0 & 0 & 0 & 0 \\
   f_2 & -1 & -1 & 0 & 1 & 0 & 0 & 0 & 0 & 0 & 0 & 0 \\
   f_3 & 0 & 0 & 0 & 0 & 0 & -1 & 0 & 0 & 0 & 0 & 0\\
   f_4 & 0 & -1 & 0 & 1 & -1 & -1 & -1 & -1 & 1 & 0 & 0\\
   f_5 & 0 & 0 & 1 & 1 & -1 & -1 & -1 & -1& 0 & 0 & -1 \\
   f_6 & 0 & 0 & 0 & 0 & 0 & 0 & 0 & 0 & 1 & 1 & 1 \\
   f_7 & 0 & 0 & 0 & 0 & 0 & 0 & 0 & 1 & -1 & -1 & 0 \\ 
  \end{array}$$
which gives equations:
$$\begin{array}{ll}
    (1) & 1=2-a_{1,2} \\
    (2) & 1=3-a_{1,2}-a_{1,4}+a_{2,4} \\
    (3) & 1=2+a_{1,5} \\
    (4) & 1=3+a_{2,4}+a_{2,5}+a_{4,5} \\
    (5) & 1=2+a_{4,5} \\
    (6) & 1=3+a_{3,4}+a_{3,5}+a_{4,5} \\
    (7) & 1=2+a_{4,5} \\
    (8) & 1=3+a_{4,5}-a_{4,7}-a_{5,7} \\
    (9) & 1=3+a_{4,6}-a_{4,7}-a_{6,7} \\
    (10)& 1=2-a_{6,7} \\
    (11)& 1=2-a_{5,6} \\
  \end{array}$$
The following matrix satisfies these equations and is positive definite:
$$\left(\begin{array}{ccccccc}
    2 & 1 & 0 & 1 & -1& 0 & 0 \\
    1 & 2 & 1 & 0 & -1& 0 & 0 \\
    0 & 1 & 2 & -1& 0 & 1 & 0 \\
    1 & 0 & -1& 2 & -1& -1& 0 \\
    -1& -1& 0 & -1& 2 & 1 & 1 \\
    0 & 0 & 1 & -1& 1 & 2 & 1\\
    0 & 0 & 0 & 0 & 1 & 1 & 2\\
  \end{array}\right)$$
And the quadratic form given by this matrix is:
\end{longversion}
$$%
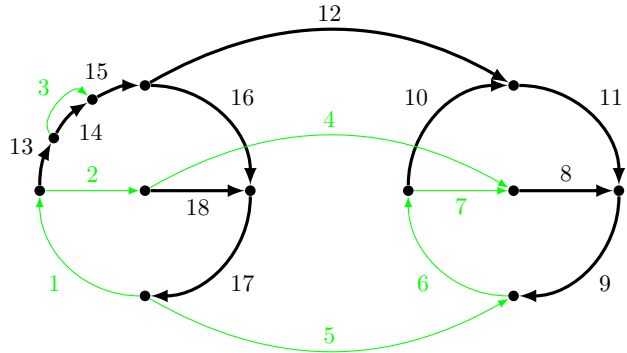
\caption{The Graph $G_{10}$}
\label{fig-g10}
\end{figure}
\end{showfigures}

\begin{longversion}
\noindent The basis for $G_{10}$ is given by:
$$\begin{array}{l}
  f_1=e_1+e_{13}+e_{14}+e_{15}+e_{16}+e_{17} \\
  f_2=e_2-e_{13}-e_{14}-e_{15}-e_{16}+e_{18} \\
  f_3=e_3-e_{14} \\
  f_4=e_4+e_8-e_{11}-e_{12}+e_{16}-e_{18} \\
  f_5=e_5-e_9-e_{11}-e_{12}+e_{16}+e_{17} \\
  f_6=e_6+e_9+e_{10}+e_{11} \\
  f_7=e_7+e_8-e_{10}-e_{11}
  \end{array}$$
Which gives us the following array:
$$\begin{array}{c|c|c|c|c|c|c|c|c|c|c|c|c|c|c|c}
        & e_8 & e_9 & e_{10} & e_{11} & e_{12} & e_{13} & e_{14} & e_{15} & e_{16} & e_{17} & e_{18} \\
\hline
   f_1 & 0 & 0 & 0 & 0 & 0 & 1 & 1 & 1 & 1 & 1 & 0 \\
   f_2 & 0 & 0 & 0 & 0 & 0 & -1 & -1 & -1 & -1 & 0 & 1 \\
   f_3 & 0 & 0 & 0 & 0 & 0 & 0 & -1 & 0 & 0 & 0 & 0 \\
   f_4 & 1 & 0 & 0 & -1 & -1 & 0 & 0 & 0 & 1 & 0 & -1 \\
   f_5 & 0 & -1 & 0 & -1 & -1 & 0 & 0 & 0 & 1 & 1 & 0 \\
   f_6 & 0 & 1 & 1 & 1 & 0 & 0 & 0 & 0 & 0 & 0 & 0 \\
   f_7 & 1 & 0 & -1 & -1 & 0 & 0 & 0 & 0 & 0 & 0 & 0 \\ 
  \end{array}$$
which gives equations:
$$\begin{array}{ll}
    (1) & 1=2+a_{4,7} \\
    (2) & 1=2-a_{5,6} \\
    (3) & 1=2-a_{6,7} \\
    (4) & 1=4+a_{4,5}-a_{4,6}+a_{4,7}-a_{5,6}+a_{5,7}-a_{6,7} \\
    (5) & 1=2+a_{4,5} \\
    (6) & 1=2-a_{1,2} \\
    (7) & 1=3-a_{1,2}-a_{1,3}+a_{2,3} \\
    (8) & 1=2-a_{1,2} \\
    (9) & 1=4-a_{1,2}+a_{1,4}+a_{1,5}-a_{2,4}-a_{2,5}+a_{4,5} \\
    (10)& 1=2+a_{1,5} \\
    (11)& 1=2-a_{2,4} \\
  \end{array}$$
The following matrix satisfies these equations and is positive definite:
$$\left(\begin{array}{ccccccc}
    2 & 1 & 1 & 1 & -1 & 0 & 0 \\
    1 & 2 & 0 & 1 & 0 & 0 & 0 \\
    1 & 0 & 2 & 0 & 0 & 0 & 0 \\
    1 & 1 & 0 & 2 & -1 & 0 & -1 \\
    -1 & 0 & 0 & -1 & 2 & 1 & 1\\
    0 & 0 & 0 & 0 & 1 & 2 & 1\\
    0 & 0 & 0 & -1 & 1 & 1 & 2\\
  \end{array}\right)$$
And the quadratic form given by this matrix is:
\end{longversion}
$$\begin{array}{rl}
  q_{G_{10}}(x_1,x_2,x_3,x_4,x_5,x_6,x_7)=& x_1^2+x_1x_2+x_1x_3+x_1x_4-x_1x_5+ \\
				 & x_2^2+x_2x_4+x_3^2+x_4^2-x_4x_5- \\
				 & x_4x_7+x_5^2+x_5x_6+x_5x_7+x_6^2+ \\
				 & x_6x_7+x_7^2 \\
  \end{array}$$

\section{Genus 8}
\label{sec:8}

As we explained in Section~\ref{subsec:genus-8}, it is sufficient to
find a $\bZ$-emm for each of the finitely many graphs obtained from
$F_{11}$--$F_{14}$ and $G_1$-$G_{10}$ by applying procedure (a) or
(b). This gives $14 \cdot \left( \binom{18}{2} + 18\right) = 2394$ graphs. 

We have written a Mathematica program for computing the $8 \times 13$
matrices for these graphs, and a 
Fortran program which uses integer arithmetics for 
finding the $\bZ$-emms. We confirmed that they exist for all of these
graphs. The lists of the matrices and the $\bZ$-emms are available at
\url{http://www.math.uga.edu/~valery/vigre2010}.

\renewcommand{\MR}[1]{}
\bibliographystyle{amsalpha}

\def\cprime{$'$}
\providecommand{\bysame}{\leavevmode\hbox to3em{\hrulefill}\thinspace}
\providecommand{\MR}{\relax\ifhmode\unskip\space\fi MR }
\providecommand{\MRhref}[2]{%
  \href{http://www.ams.org/mathscinet-getitem?mr=#1}{#2}
}
\providecommand{\href}[2]{#2}

\end{document}